\documentclass[a4paper,12pt]{article}
\usepackage{amsmath}
\usepackage{amsthm}
\usepackage{amsfonts}
\usepackage{amssymb}
\usepackage{amstext}
\usepackage{amsopn}
\usepackage{bbm}
\usepackage{graphicx}
\usepackage{mathrsfs}
\usepackage{multirow}

\usepackage{color}
\usepackage{mathtools}

\newcommand{\al}{\alpha}		%alpha
\newcommand{\vep}{\varepsilon}	%variant epsilon
\newcommand{\ze}{\zeta}		%zeta
\newcommand{\lam}{\lambda}	%lambda
\newcommand{\sig}{\sigma}	%sigma
\newcommand{\ome}{\omega}	%omega
\newcommand{\Ome}{\Omega}	%Omega

\newcommand{\bIn}{\mathbbm 1}		%indicator function
\newcommand{\bC}{\mathbb C}		%complex numbers
\newcommand{\bD}{\mathbb D}		%unit disk
\newcommand{\bE}{\mathbb E}		%expectation
\newcommand{\bM}{\mathbb M}		%space of probability measures
\newcommand{\bN}{\mathbb N}		%natural numbers
\newcommand{\bP}{\mathbb P}		%probability
\newcommand{\bR}{\mathbb R}		%real numbers
\newcommand{\cB}{\mathcal B}		%Borel sigma-algebra
\newcommand{\cF}{\mathcal F}		%sigma-algebra
\newcommand{\cL}{\mathcal L}		%law / Laplace transform
\newcommand{\UC}{\partial\bD}		%unit circle
\newcommand{\re}{{\rm Re}}		%real part
\newcommand{\im}{{\rm Im}}		%imaginary part
\newcommand{\diag}{{\rm diag}}		%diagonal matrix
\newcommand{\tr}{{\rm tr}}		%matrix trace
\newcommand{\abs}[1]{\left|{#1}\right|}			%absolute value sign
\newcommand{\rob}[1]{\left({#1}\right)}		%round brackets
\newcommand{\sqb}[1]{\left[{#1}\right]}		%square brackets
\newcommand{\cub}[1]{\left\{{#1}\right\}}		%curly brackets
\newcommand{\inp}[2]{\left\langle#1,#2\right\rangle}		%inner product
\newcommand{\floor}[1]{\left\lfloor{#1}\right\rfloor}	%floor function
\newcommand{\fr}{\frac}		%fraction
\newcommand{\dfr}{\dfrac}		%fraction (display style)
\newcommand{\nth}[1]{\fr{1}{#1}}			%1/n
\newcommand{\half}{\fr{1}{2}}			%1/2
\newcommand{\inclusion}{\hookrightarrow}		%inclusion map
\newcommand{\wto}{\rightharpoonup}		%weak convergence
\newcommand{\arr}{\rightarrowtail}		%convergence to a set
\newcommand{\To}{\Longrightarrow}		%only if
\newcommand{\antti}{\text{ \ as \ }n\to\infty}	%as n tends to infinity
\newcommand{\ntti}{n\to\infty}	%n tends to infinity
\newcommand{\wt}{\widetilde}				%wide tilde
\newcommand{\ol}{\overline}				%overline
\newcommand{\ub}{\underbrace}				%underbrace
\newcommand{\dst}{\displaystyle}			%displaystyle
\definecolor{w}{rgb}{1,1,1}					%white color
\definecolor{r}{rgb}{0,0,0}					%red color
\definecolor{g}{rgb}{0,0,0}					%green color
\definecolor{b}{rgb}{0,0,0}					%blue color
\theoremstyle{definition}
	\newtheorem{dfn}{Definition}
	
	\newtheorem*{remark}{Remark}
\theoremstyle{plain}
	\newtheorem{lem}[dfn]{Lemma}
	\newtheorem{pro}[dfn]{Proposition}
	\newtheorem{thm}[dfn]{Theorem}
	\newtheorem{cor}[dfn]{Corollary}
	\newtheorem{cnj}{}

\begin{document}

\title{\bf Higher Order, Polar and Sz.-Nagy's Generalized Derivatives of Random Polynomials with Independent and Identically Distributed Zeros\\on the Unit Circle}
\author{\textsc{
\begin{tabular}{c}
Pak-Leong Cheung${}^{1}$\thanks{
{\color{b}Partially supported by a graduate studentship of HKU and the RGC grants HKU 706411P and HKU 703313P.}
\mbox{\hspace{4truecm}}}
, Tuen Wai Ng${}^{1}$\thanks{
{\color{b}Partially supported by the RGC grant HKU 704611P and 703313P.}
\mbox{\hspace{0.9truecm}}}
,\\Jonathan Tsai${}^{1}$\thanks{
{\color{b}Partially supported by a faculty postdoctoral fellowship of the Faculty of Science, HKU.}
\mbox{\hspace{0.9truecm}}}
\ and S.C.P. Yam${}^{2}$\thanks{
{\color{b}Acknowledges The Chinese University of Hong Kong Direct Grant 2011/2012 Project ID: 2060444.}
\mbox{\hspace{0.9truecm}}}
\end{tabular}}}
 
\date{}
\maketitle

\centerline{\textbf{\today}}

\begin{figure}[b]
\rule[-2.5truemm]{5cm}{0.1truemm}\\[2mm]
{\footnotesize  
{\color{b}2000 {\it Mathematics Subject Classification: Primary} 30C10, 30B20, 30C15, {\it Secondary} 60B10, 60G57.}
\par {\it Key words and phrases.} 
{\color{b}Random polynomial, zero distribution, polar derivative, Sz.-Nagy's generalized derivative, random measure}
\par\noindent 1.
Department of Mathematics,
The University of Hong Kong,
Pokfulam, Hong Kong.\\
2.
Department of Statistics,
The Chinese University of Hong Kong,
Shatin, Hong Kong.
\smallskip
\par\noindent E-mail: {\tt mathcpl@connect.hku.hk}, {\tt ntw@maths.hku.hk}, {\tt jonathan.tsai@rocketmail.com}, {\tt scpyam@sta.cuhk.edu.hk}}
\end{figure}

\begin{quotation}
\noindent{\textbf{Abstract}. For random polynomials with i.i.d. (independent and identically distribu-ted) zeros following any common probability distribution $\mu$ with support contained in the unit circle, the empirical measures of the zeros of their first and higher order derivatives will be proved to converge weakly to $\mu$ a.s. (almost sure(ly)). This, in particular, completes a recent work of Subramanian on the first order derivative case where $\mu$ was assumed to be non-uniform. The same a.s. weak convergence will also be shown for polar and Sz.-Nagy's generalized derivatives, on some mild conditions.}
\end{quotation}

% =====================================================================
\bigskip
\section{Introduction}\label{sec1}%sec1

\bigskip
The study of zero distribution of random polynomials {\color{b}(and random power series)} has a long history and is currently a very active research area (see the references in \cite{HoKrPerVi}, \cite{HuNi}, \cite{PerVi}, \cite{ShepVa} and \cite{ShZe}). Traditionally, the randomness in these polynomials {\color{b}comes from the probability distribution followed by their coefficients, i.e.
$$A_0+A_1z+\cdots+A_nz^n,$$
where $A_0,A_1,\dots,A_n$ are complex-valued random variables (for instance, i.i.d. Gaussian).

\bigskip
Instead of the polynomials' coefficients, one may opt to introduce randomness in their zeros}, and then investigate the locations of their critical points (relative to these zeros). Such a study was initiated by Rivin and the late Schramm in 2001, but only until 2011, Pemantle and Rivin \cite{PemRi} proposed a precise probabilistic framework {\color{b}(see Table \ref{tab2.101}).

\bigskip
\begin{table}[h]%---------- TABLE 1 ----------
\centering
\begin{tabular}{ccc}
\bf Random polynomials	&\bf Traditional	&\bf Pemantle--Rivin\\
\hline
Prescribe randomness to	&Coefficients	&Zeros\\
\hline
To study				&Zeros		&Critical points
\end{tabular}
\caption{Summary of the traditional and the Pemantle--Rivin frameworks of random polynomials\label{tab2.101}}%{tab2.101}
\end{table}

\bigskip\noindent
Their framework is a case of the following setting:}

\bigskip
Let $\mu_S$ be a probability measure on $\bC$ (with the Borel $\sig$-algebra $\cB$) {\it supported in} a closed set $S\subset\bC$ (meaning that $\text{supp}\,\mu_S\subset S$ or, equivalently, $\mu_S(S)=1$) and
\begin{equation}\label{eq1}%{eq1}
Z=Z_1,Z_2,\dots\phantom{}\sim\mu_S
\end{equation}
be i.i.d. complex-valued random variables on a probability space $(\Ome,\bP)$. Fix any $k\in\bN$. For each $\ome\in\Ome$ and $n\geq k+1$, construct the polynomial
\begin{equation}\label{eq2}%{eq2}
P_n(\ome)(z):=(z-Z_1(\ome))\cdots(z-Z_n(\ome)).
\end{equation}
By {\color{b}regarding
$$P_n(z)=(z-Z_1)\cdots(z-Z_n)$$
as} a random variable on $\Ome$ into the space of polynomials, and taking \eqref{eq1} into account so that all $Z_j\in S$ a.s., we call $P_n$ a {\it random polynomial with i.i.d. zeros $Z_1,\dots,Z_n$ on $S$}.

\bigskip
{\color{b}Associate $n-k$ points
\begin{equation}\label{eq3}%{eq3}
W_{n,1}(\ome),W_{n,2}(\ome),\dots,W_{n,n-k}(\ome)\in\bC,
\end{equation}
to the polynomial $P_n(\ome)$ in \eqref{eq2}. In the original Pemantle--Rivin framework, \eqref{eq3} are taken as the critical points
$$W^{(1)}_{n,1}(\ome),W^{(1)}_{n,2}(\ome),\dots,W^{(1)}_{n,n-1}(\ome)$$
of $P_n(\ome)$ ($k=1$). Later in this chapter, \eqref{eq3} will also be taken as the zeros
$$W^{(k)}_{n,1}(\ome),W^{(k)}_{n,2}(\ome),\dots,W^{(k)}_{n,n-k}(\ome)$$
of higher order derivatives ($k\geq2$), those of polar derivatives ($k=1$) and Sz.-Nagy's generalized derivatives ($k=1$) of $P_n(\ome)$ (to be made precise in the paragraph following Corollary \ref{cor4}). Define} the {\it empirical measure} of these $n-k$ points {\color{b}in \eqref{eq3}} by
\begin{equation}\label{eq4}%{eq4}
\mu_{S,n}(\ome):=\nth{n-k}\sum_{j=1}^{n-k}\delta_{W_{n,j}(\ome)},
\end{equation}
where $\delta_x$ is the Dirac measure with support at $\cub{x}$. Note that $\mu_{S,n}(\ome)$ is also a probability measure on $(\bC,\cB)$, so we may discuss the weak convergence
\begin{equation}\label{eq5}%{eq5}
\mu_{S,n}(\ome)\wto\mu_S\antti
\end{equation}
of probability measures. Moreover, endow the space $\bM$ of all probability measures on $(\bC,\cB)$ with the metric topology of weak convergence. Then, regarding
$$\mu_{S,n}:=\nth{n-k}\sum_{j=1}^{n-k}\delta_{W_{n,j}}:(\Ome,\bP)\to\bM$$
and $\mu_S$ as random variables (the latter being a constant function) on $\Ome$ into $\bM$, we may discuss various modes (e.g. almost sure, in probability, etc.) of the (weak) convergence of {\it random (empirical) measures}
\begin{equation}\label{eq6}%{eq6}
\mu_{S,n}\wto\mu_S\antti
\end{equation}
{\color{b}Two of these modes will concern us:

\bigskip
\begin{dfn}\label{dfn2.102}%{dfn2.102}
By that the {\it weak convergence \eqref{eq6} holds almost surely}, we mean, as usual, that
\begin{equation}\label{eq101}%{eq101}
\bP\cub{\ome\in\Ome:\mu_{S,n}(\ome)\wto\mu_S\text{ as }\ntti}=1
\end{equation}
or, equivalently,
\begin{equation}\label{eq102}%{eq102}
\bP\cub{\ome\in\Ome:\pi(\mu_{S,n}(\ome),\mu_S)\to0\text{ as }\ntti}=1.
\end{equation}
And by that the {\it weak convergence \eqref{eq6} holds in probability}, we mean that: For any $\vep>0$,
\begin{equation}\label{eq103}%{eq103}
\bP\cub{\ome\in\Ome:\pi\rob{\mu_{S,n}(\ome),\mu_S}<\vep}\to1\antti.
\end{equation}
In \eqref{eq102} and \eqref{eq103}, $\pi$ denotes the Prohorov metric
\begin{equation}\label{eq38}%{eq38}
\begin{array}{rl}
\pi(m',m''):=\inf\{\vep>0:&m'(A)\leq m''(A^\vep)+\vep\text{ \ and}\\
&m''(A)\leq m'(A^\vep)+\vep\text{ \ for all \ }A\in\cB\}
\end{array}
\end{equation}
$$\text{for \ }m',m''\in\bM,$$
where \smash{$A^\vep:=\bigcup_{a\in A}\bD_\vep(a)$}, which gives the metric topology on $\bM$ (\cite[Theorem 6.8, p.73]{Bi}) just mentioned.
\end{dfn}

\bigskip
\begin{remark}
Note that \eqref{eq102} implies \eqref{eq103}. And, convergence in probability only guarantees almost sure convergence for some subsequence.
\end{remark}

\bigskip\noindent
Throughout} this article, all `a.s.' and `in probability' statements and arguments are said with respect to $(\Ome,\bP)$, and the zeros of any holomorphic function (thus any polynomial) are counted with multiplicities.

\bigskip
When $W_{n,1}(\ome),\dots,W_{n,n-1}(\ome)$ in \eqref{eq3} (with $k=1$) are the critical points \smash{\color{b}$W^{(1)}_{n,1}(\ome),$} \smash{\color{b}$\dots,W^{(1)}_{n,n-1}(\ome)$} of $P_n(\ome)$ in \eqref{eq2}, {\color{b}we rewrite \eqref{eq4} as
$$\mu^{(1)}_{S,n}(\ome)=\nth{n-1}\sum_{j=1}^{n-1}\delta_{W^{(1)}_{n,j}}(\ome).$$
In such case, the scenario as set in Section \ref{sec1}} reduces to Pemantle and Rivin's study in \cite{PemRi}. They conjectured{\color{b}, in view of a couple of supportive examples, that:

\bigskip
\begin{cnj}%Conjecture A
\hspace{-.66cm}\colorbox{w}{\bf Conjecture A.} {\bf(Pemantle--Rivin conjecture)} {\rm(\cite[Conjecture 1.1]{PemRi})} The weak convergence \eqref{eq6} {\color{b}with \smash{$\mu_{S,n}=\mu^{(1)}_{S,n}$}} holds in probability.
\end{cnj}

\bigskip\noindent
In} fact, they proved their conjecture when $\mu_S$ has finite 1-energy and if
$$\mu_S\cub{z\in\bC:\int_{\ze\in\bC}\fr{d\mu_S(\ze)}{z-\ze}=0}=0,$$
and also pointed out that their method of proof does not apply to the case when $S=\UC$ with uniform $\mu_{\UC}$, which has infinite 1-energy. Yet, in such case, they deduced that \smash{\color{b}$\mu^{(1)}_{\UC,n}$} converges to $\UC$ (which intuitively means that `\smash{\color{b}$\mu^{(1)}_{\UC,n}$} tends to be supported in $\UC$ when $n$ gets large', and the precise definition will be stated in Section \ref{sec2}) in probability (necessary for the conjecture to be true; see Proposition \ref{pro2} below) from a result of Peres and Vir\'ag \cite{PerVi} on determinantal point process. Soon later, Subramanian \cite{Sub} followed up on \cite{PemRi} and showed that {\color{b}Conjecture A} is true when $S=\UC$. In the same period, Kabluchko \cite{Ka} actually confirmed the conjecture in full generality, i.e. for $S=\bC$, using a natural connection with logarithmic potential theory.

\bigskip
A harder question is the almost sureness (instead of `in probability') of {\color{b}\eqref{eq6}:

\bigskip
\begin{cnj}%Conjecture B
\hspace{-.66cm}\colorbox{w}{\bf Conjecture B.} {\bf(Strong Pemantle--Rivin conjecture)} The weak convergence \eqref{eq6} {\color{b}with \smash{$\mu_{S,n}=\mu^{(1)}_{S,n}$}} holds almost surely.
\end{cnj}

\bigskip\noindent
Subramanian} \cite{Sub} proved it for the case that $S=\UC$ with non-uniform $\mu_{\UC}$ (implicit in \cite[Theorem 2.3]{Sub}). The methods in its proof and the remaining case of uniform $\mu_{\UC}$ have motivated the present work. Besides, rather than only considering the critical points in {\color{b}Conjecture B}, we extend the problem to the zeros of higher order, polar and Sz.-Nagy's generalized derivatives of $P_n(\ome)$ in \eqref{eq2}. Precisely, we will prove that
$$\mu_{\UC,n}\wto\mu_{\UC}\antti\text{ \ a.s.}$$
(i.e. the almost sureness of the weak convergence \eqref{eq6} with $S=\UC$) for the case when $W_{n,1}(\ome),\dots,W_{n,n-k}(\ome)$ in \eqref{eq3} are the zeros of
\begin{enumerate}
\item[(i)]	the {\color{b}$k$-th} order derivative of $P_n(\ome)$ with non-uniform $\mu_{\UC}$ and $k>1$ (Corollary \ref{cor4}).
\item[(ii)]	the polar derivative of $P_n(\ome)$ with respect to any {\color{b}$\xi\in\bC$ with $\abs{\xi}\gg1$} {\color{b}(more precisely, $\abs{\xi}>1$ if $\mu_{\UC}$ is uniform and $\abs{\xi}>\max\cub{1,1/\abs{\bE[Z]}}$ if $\mu_{\UC}$ is non-uniform)} (Corollary \ref{cor5}).
\item[(iii)]	any Sz.-Nagy's generalized derivative of $P_n(\ome)$ given by degenerate (i.e. a.s. constant) random coefficients which are uniformly bounded a.s. and with non-uniform $\mu_{\UC}$ (Corollary \ref{cor5}).
\item[(iv)]	the {\color{b}$k$-th} order derivative of $P_n(\ome)$ with uniform $\mu_{\UC}$ and $k\geq1$ (Corollary \ref{cor6}).
\end{enumerate}

\bigskip\noindent
{\color{b}We shall rewrite \eqref{eq6} as
$$\mu^{(k)}_{\UC,n}\wto\mu_{\UC}\antti$$
for cases (i) and (iv), and as
$$\mu^\xi_{\UC,n}\wto\mu_{\UC}\antti$$
for case (ii).} In particular, (i) and (iv) combine with the aforesaid result implicit in \cite[Theorem 2.3]{Sub} to form the {\color{b}collective result}:

%---------- THEOREM 1 ----------
\bigskip
\begin{thm}\label{thm1}
Fix any $k\in\bN$. Let $W_{n,1}(\ome),\dots,W_{n,n-k}(\ome)$ in \eqref{eq3} be the zeros of the {\color{b}$k$-th} order derivative of $P_n(\ome)$ in \eqref{eq2}. Then, \smash{\color{b}$\mu^{(k)}_{\UC,n}\wto\mu_{\UC}$} as $\ntti$ a.s..
\end{thm}

\bigskip\noindent
which completes the discussion on {\color{b}Conjecture B} when $S=\UC$ in \cite{Sub}, and generalizes it to higher order derivatives (summarized in Table \ref{tab1} below). Note that this result is non-trivial --- it cannot be obtained, for instance, by repeated applications of the results in \cite{Sub} because the critical points of $P_1,P_2,\dots$ in \eqref{eq2} derived from $Z_1,Z_2,\dots$ in \eqref{eq1} may not be i.i.d..

%---------- TABLE 1 ----------
\bigskip
\begin{table}[h]
\centering
\begin{tabular}{cccc}
$\mu_{\UC}$&\!\!\!\!\!&$k=1$&$k>1$\\
\cline{1-1}\cline{3-4}
non-uniform&\!\!\!\!\!&\cite[Theorem 2.3]{Sub}&Corollary \ref{cor4}\\
\cline{1-1}\cline{3-4}
uniform&\!\!\!\!\!&\multicolumn{2}{c}{Corollary \ref{cor6}}
\end{tabular}
\caption{{\color{b}Conjecture B} when $S=\UC$ and generalization to higher order derivatives\label{tab1}}%{tab1}
\end{table}

\bigskip
All the results (i)--(iv) will be consequences of Proposition \ref{pro3}, which, under the assumption that in \eqref{eq3} all {\color{b}$W_{n,j}\in\ol{\bD_M}$} a.s. {\color{b}for some constant $M>0$}, captures the essentials of the a.s. weak convergence \eqref{eq6} for the natural case when $S=\UC$. This proposition can be regarded as originating from \cite[Theorem 2.3]{Sub}. Despite this connection, the present result differs from the earlier by, on one hand, its higher applicability (for instances, to the zeros of higher order and other types of derivatives) and by, on the other hand, featuring L\'evy's continuity theorem (as far as weak convergence of probability measures is concerned). Via this remarkable theorem {\color{b}of L\'evy}, Section \ref{sec3.1} will set up a straightforward formulation of our problem in the general setting, so that the weak convergence \eqref{eq5} can then be discussed by examining characteristic functions expressed explicitly in terms of $W_{n,j}$ in \eqref{eq3}.

\bigskip
In Corollaries \ref{cor4} and \ref{cor5}, we will verify that when $W_{n,j}$ in \eqref{eq3} are respectively the zeros of (i), (ii) and (iii) above, the assumption and sufficient condition in Proposition \ref{pro3} are satisfied. These verifications will be done by elaborating the following two main ideas in \cite{Sub}:
\begin{enumerate}
\item[$\cdot$]	{\it Construct a non-identically zero holomorphic function on $\bD$ whose number of zeros lying in any compact subset of $\bD$ can be related to that of the derivative in (i), (ii) and (iii) respectively by Hurwitz's theorem;}
\item[$\cdot$]	{\it The asymptotic relationship between the means of the powers of the zeros and the critical points of a polynomial when its degree tends to infinity (implied by Lemma \ref{lem7} below).}
\end{enumerate}
However, the same construction of holomorphic function in (iv) (now $\mu_{\UC}$ is uniform instead of being non-uniform in (i)) would result in the identically zero function, to which Hurwitz's theorem is not applicable. In this exceptional case, we discovered (from the proof of (ii)) a uniquely nice behaviour of the zeros of polar derivative. This property can be passed to those of the {\color{b}$k$-th} order derivative via approximating ordinary derivatives by polar derivatives, thus proving Corollary \ref{cor6}:

%---------- TABLE 1 ----------
\bigskip
\begin{table}[h]
\centering
{\color{b}\begin{tabular}{ccccc}
&\multicolumn{4}{c}{$\mu_{\UC,n}\wto\mu_{\UC}$}\\
\cline{3-5}
$\mu_{\UC}$&\!\!\!\!\!&$\mu^{(1)}_{\UC,n}$&$\mu^{(k)}_{\UC,n}$, $k>1$&$\mu^\xi_{\UC,n}$, $\abs{\xi}\gg1$\\
\cline{1-1}\cline{3-5}
non-uniform&\!\!\!\!\!&\cite[Theorem 2.3]{Sub}&Corollary \ref{cor4}&\multirow{2}{*}{Corollary \ref{cor5}}\\
\cline{1-1}\cline{3-4}
uniform&\!\!\!\!\!&\multicolumn{2}{c}{\phantom{$\Longleftarrow$}\qquad\qquad Corollary \ref{cor6}\qquad\qquad$\Longleftarrow$}
\end{tabular}}
\caption{Towards $\mu^{(k)}_{\UC,n}\wto\mu_{\UC}$ with uniform $\mu_{\UC}$\label{tab2.102}}%{tab2.102}
\end{table}

\bigskip
After presenting the above main contents by the end of Section \ref{sec3}, we remark on a rather rudimentary connection between the weak convergence \eqref{eq6} and a classical question of the locations of critical points relative to zeros of polynomials (Corollary \ref{cor8}, and a similar remark has also been mentioned in \cite{Ka}) in Section \ref{sec4}.

% ===========================================================================
\bigskip
\section{Main result and corollaries}\label{sec2}

\bigskip
First of all, we mention a mode of convergence closely related to the weak convergence of probability measures. We explicitly define this notion which has already appeared in \cite{PemRi} and \cite{Sub} for readers' easy reference. Let $S\subset\bC$ be a closed set. Then, a sequence $m_1,m_2,\dots$ of probability measures on $(\bC,\cB)$ is said to {\it converge to $S$} (denoted as $m_n\arr S$) as $\ntti$ {\color{b}if
\begin{equation}\label{convtoset}%{convtoset}
m_n(S^c)\to0\antti.
\end{equation}
It} is immediate from definition and the portmanteau theorem (\cite[Theorem 2.1]{Bi}, \cite[Theorem 11.1.1]{Dud}, \cite[Theorem 3.2.5]{Dur}) to have:

%---------- PROPOSITION 2 ----------
\bigskip
\begin{pro}\label{pro2}
Let $m,m_1,m_2,\dots$ be probability measures on $(\bC,\cB)$. If $m$ is supported in a closed set $S\subset\bC$, then
$$m_n\wto m\antti\qquad\To\qquad m_n\arr S\antti.$$
\end{pro}

%---------- PROOF OF PROPOSITION 2 ----------
\bigskip
\begin{proof}
{\color{b}Since $m_n\wto m$ as $\ntti$, it follows from the portmanteau theorem that
$$\limsup_{\ntti}m_n(B)\leq m(B)\text{ \ for any closed subset \ }B\subset\bC.$$
Exhaust the open set $S^c$ by (an increasing sequence of) compact subsets $B_k$ (i.e. \smash{$B_1\Subset B_2\Subset\cdots\subset S^c=\bigcup_{k=1}^\infty B_k$}), then we actually have
$$\limsup_{\ntti}m_n(B_k)\leq m(B_k)=0\text{ \ for all \ }k$$
because $\text{supp}\,m\subset S$. Let \smash{$\limsup_{\ntti}m_n(S^c)$} be realized by subsequence\newline\smash{$\cub{m_{n_i}(S^c)}_{i=1}^\infty$}. Then, we have}
\begin{align*}
\limsup_{\ntti}m_n(S^c)&=\lim_{i\to\infty}m_{n_i}(S^c)=\lim_{i\to\infty}\lim_{k\to\infty}m_{n_i}(B_k)\\
&=\lim_{k\to\infty}\lim_{i\to\infty}m_{n_i}(B_k)\leq\lim_{k\to\infty}\limsup_{\ntti}m_n(B_k)\leq0.\qedhere
\end{align*}
\end{proof}

\bigskip
\begin{remark}
{\color{b}Note that \eqref{convtoset} is equivalent to `$m_n(B)\to0$ as $\ntti$ for any closed subset $B\subset S^c$' via a compact exhaustion of the open set $S^c$ as in the above proof: $\limsup_nm_n(S^c)=\lim_im_{n_i}(S^c)=\lim_i\lim_km_{n_i}(B_k)=\lim_k\lim_im_{n_i}(B_k)=0$.}
\end{remark}

\bigskip
From now on (except in Section \ref{sec3.1}), $\mu_S$ is assumed to be supported in $\UC$, i.e. $S=\UC$ and all $\abs{Z_j}=1$ a.s.. For brevity, we suppress all subscripts `$\UC$' by {\color{b}writing
$$\mu=\mu_{\UC},\ \ \mu_n=\mu_{\UC,n},\ \ \mu^{(k)}_n=\mu^{(k)}_{\UC,n}\text{ \ and \ }\mu^\xi_n=\mu^\xi_{\UC,n}$$
in} what follows. Also, we use the usual notation $\bE_m$ for the expectation with respect to a probability measure $m$, with $\bE_\bP$ abbreviated as $\bE$. According to Proposition \ref{pro2}, $\mu_n(\ome)\arr\UC$ is necessary for $\mu_n(\ome)\wto\mu$ in \eqref{eq5} when $\ntti$. Indeed, we also have

%---------- PROPOSITION 3 ----------
\bigskip
\begin{pro}\label{pro3}
Let $\mu$ be supported in $\UC$. Assume that {\color{b}there is a constant $M>0$ such that in \eqref{eq3}
\begin{equation}\label{eq7}%{eq7}
\abs{W_{n,j}}\leq M\ \ \forall\ n,j\text{ \ a.s..}
\end{equation}
Then}, $\mu_n\wto\mu$ as $\ntti$ a.s. in \eqref{eq5} if and only if the following two conditions hold:
\begin{center}
\begin{tabular}{rl}
(i)&$\mu_n\arr\UC$ \ as \ $\ntti$ \ a.s.\vspace{.1cm}\\
(ii)&$\dst\nth{n-k}\sum_{j=1}^{n-k}{W_{n,j}}^p\to\bE\sqb{Z^p}$ \ as \ $\ntti\ \ \forall\ p\in\bN$ \ a.s.
\end{tabular}
\end{center}
\end{pro}

\bigskip
Making use of a key lemma, namely Lemma \ref{lem7} (see \cite[Proposition 3.2]{Sub}) which is to be proved in Section \ref{sec3.3}, Proposition \ref{pro3} leads to the following three corollaries:

%---------- COROLLARY 4 ----------
\bigskip
\begin{cor}\label{cor4}
Let $\mu$ be supported in $\UC$ and $W_{n,1}(\ome),\dots,W_{n,n-k}(\ome)$ in \eqref{eq3} (with $k>1$) be the zeros of the {\color{b}$k$-th} order derivative of $P_n(\ome)$ in \eqref{eq2}. Then, \smash{\color{b}$\mu^{(k)}_n\wto\mu$} as $\ntti$ a.s. if $\mu$ is non-uniform.
\end{cor}

\bigskip\noindent
Before stating the next corollary, we need to define the following two types of derivatives for a polynomial
\begin{equation}\label{eq8}%{eq8}
P(z):=(z-z_1)\cdots(z-z_n).
\end{equation}
Let $Q$ be a polynomial given by
\begin{equation}\label{eq9}%{eq9}
Q(z):=P(z)\sum_{j=1}^n\fr{\lam_j}{z-z_j}=\rob{\sum_{j=1}^n\lam_j}z^{n-1}+\cdots,
\end{equation}
where $\lam_1,\dots,\lam_n\in\bC$. Assume $\deg Q=n-1$ so that $\sum_{j=1}^n\lam_j\neq0$, then $Q$ is called
\begin{enumerate}
\item[$\cdot$]	the {\it polar derivative} $D_\xi P$ of $P$ {\it with respect to $\xi\in\bC$} if
$$\lam_j=\xi-z_j$$
(\cite{CheuNg2}, \cite[p.184]{MilMitRas}, \cite[p.97]{RahSc}, \cite[p.185]{Shei}). In such case,
$$Q(z)=D_\xi P(z)=nP(z)-(z-\xi)P'(z)=\rob{n\xi-\sum_{j=1}^nz_j}z^{n-1}+\cdots.$$
\item[$\cdot$]	a {\it Sz.-Nagy's generalized derivative} of $P$ if 
$$\sum_{j=1}^n\lam_j=n\text{ \ and \ }\lam_j>0$$
(\cite{CheuNg2}, \cite[p.115]{RahSc}, Sz.-Nagy's original paper \cite{Sz}). In particular, when $\lam_j=1$, we have $Q=P'$.
\end{enumerate}

%---------- COROLLARY 5 ----------
\bigskip
\begin{cor}\label{cor5}
Let $\mu$ be supported in $\UC$ and $W_{n,1}(\ome),\dots,W_{n,n-1}(\ome)$ in \eqref{eq3} (with $k=1$) be the zeros of the polynomial
$$Q_n(\ome)(z):=P_n(\ome)(z)\sum_{j=1}^n\fr{\lam_{n,j}(\ome)}{z-Z_j(\ome)},$$
where $P_n(\ome)$ is as in \eqref{eq2} and $\lam_{n,j}:\Ome\to\bC$ ($n=2,3,\dots$, $j=1,\dots,n$) are such that
$$\sum_{j=1}^{n}\lam_{n,j}\neq0\ \ \forall\ n\text{ \ a.s.}$$
so that $\deg Q_n=n-1$ a.s.. Assume that $W_{n,j}$ satisfy \eqref{eq7}. Then, $\mu_n\wto\mu$ as $\ntti$ a.s. if the following three conditions hold:
\begin{center}
\begin{tabular}{rl}
(i)&$\dst\nth{n}\sum_{j=1}^{n}\ol{\lam_{n,j}}{Z_j}^{m+1}\to b_m$ \ as \ $\ntti$ \ a.s.\vspace{.1cm}\\
(ii)&$\dst\nth{n}\sum_{j=1}^{n}\abs{\lam_{n,j}}\leq K\ \ \forall\ n$ \ a.s.\vspace{.1cm}\\
(iii)&$\dst\nth{n}\abs{\sum_{j=1}^{n}\lam_{n,j}}\geq\delta\ \ \forall\ n$ \ a.s.
\end{tabular}
\end{center}
where $K,\delta>0$ and $b_0,b_1,\dots\phantom{}\in\bC$ are constants, and at least one of which is non-zero. This result is applicable when $Q_n$ is
\begin{enumerate}
\item[$\cdot$]	the polar derivative of $P_n$ with respect to any {\color{b}$\xi\in\bC$ with $\abs{\xi}\gg1$} {\color{b}(more precisely, $\abs{\xi}>1$ if $\mu_{\UC}$ is uniform and $\abs{\xi}>\max\cub{1,1/\abs{\bE[Z]}}$ if $\mu_{\UC}$ is non-uniform)}.
\item[$\cdot$]	any Sz.-Nagy's generalized derivative of $P_n$ such that each $\lam_{n,j}$ is degenerate (i.e. a.s. constant; but the $\lam_{n,j}$'s need not take the same constant value a.s. for different $n$ or $j$) and there is a constant $M>0$ such that all $\lam_{n,j}\leq M$ a.s. and if $\mu$ is non-uniform.
\end{enumerate}
\end{cor}

%---------- COROLLARY 6 ----------
\bigskip
\begin{cor}\label{cor6}
Let $\mu$ be supported in $\UC$ and $W_{n,1}(\ome),\dots,W_{n,n-k}(\ome)$ in \eqref{eq3} be the zeros of the {\color{b}$k$-th} order derivative of $P_n(\ome)$ in \eqref{eq2}. Then, \smash{\color{b}$\mu^{(k)}_n\wto\mu$} as $\ntti$ a.s. if $\mu$ is uniform.
\end{cor}

%%%%%%%%%% 3. PROOFS %%%%%%%%%%
\bigskip
\section{Proofs}
\label{sec3}

%%%%% 3.1 FORMULATION USING LEVY'S CONTINUITY THEOREM %%%%%
\bigskip
\subsection{Formulation using L\'evy's continuity theorem}
\label{sec3.1}

\bigskip
In this subsection, we establish a natural connection between the question of the weak convergence
$$\mu_{S,n}(\ome)\wto\mu_S\antti$$
in \eqref{eq5} and L\'evy's continuity theorem (\cite[Theorem 9.8.2]{Dud}, \cite[Theorem 3.9.4]{Dur}), which states that:
\begin{quote}
\it Let $X,X_1,X_2,\dots$ be $\bR^d$-valued random variables. Then, $\cL_{X_n}\wto\cL_X$ as $\ntti$ if and only if
$$\varphi_{X_n}(t)\to\varphi_X(t)\antti\ \ \forall\ t\in\bR^d,$$
where $\cL_X,\cL_{X_1},\cL_{X_2},\dots$ and $\varphi_X,\varphi_{X_1},\varphi_{X_2},\dots$ are the laws and the characteristic functions of $X,X_1,X_2,\dots$ respectively.
\end{quote}
For each $\ome\in\Ome$ and $n\geq k+1$, define probability space $(\Ome_n(\ome),\cF_n(\ome),$ $\bP_n(\ome))$ by
$$\Ome_n(\ome):=\cub{W_{n,1}(\ome),\dots,W_{n,n-k}(\ome)}\text{ \ (as in \eqref{eq3})},$$
\begin{equation}\label{eq10}%{eq10}
\cF_n(\ome):=2^{\Ome_n(\ome)}\text{ \ and \ }\bP_n(\ome)(A):=\fr{\abs{A}}{n-k}\text{ \ for \ }A\in\cF_n(\ome).
\end{equation}
Consider the inclusion map
\begin{equation}\label{eq11}%{eq11}
W_n(\ome):\Ome_n(\ome)\inclusion\bC,\ W_{n,j}(\ome)\mapsto W_{n,j}(\ome)\ (j=1,\dots,n-k).
\end{equation}
Note that for any {\color{b}$B\in\cB(\bC)$},
\begin{align*}
&\bP_n(\ome)\cub{\ome'\in\Ome_n(\ome):W_n(\ome)(\ome')\in B}=\bP_n(\ome)\cub{\ome'\in\Ome_n(\ome):\ome'\in B}\\
&\qquad\qquad\qquad\quad=\nth{n-k}\sum_{j=1}^{n-k}\bIn_B(W_{n,j}(\ome))=\nth{n-k}\sum_{j=1}^{n-k}\delta_{W_{n,j}(\ome)}(B)=\mu_{S,n}(\ome)(B),
\end{align*}
thus
\begin{equation}\label{eq12}%{eq12}
\mu_{S,n}(\ome)=\cL_{W_n(\ome)}
\end{equation}
which is the law of $W_n(\ome)$ in \eqref{eq11}. Therefore, it follows from L\'evy's continuity theorem that the weak convergence \eqref{eq5} is equivalent to
\begin{equation}\label{eq13}%{eq13}
\varphi_{W_n(\ome)}(t)\to\varphi_Z(t)\antti\ \ \forall\ t\in\bC.
\end{equation}
Indeed, by \eqref{eq12} and \eqref{eq4}, we can write
\begin{align}
\varphi_{W_n(\ome)}(t)&:=\bE_{\bP_n(\ome)}\!\sqb{e^{i\inp{t}{W_n(\ome)}}}=\int_{\ze\in\bC}e^{i\inp{t}{\ze}}\,d\mu_{S,n}(\ome)(\ze)\nonumber\\
&\phantom{:}=\nth{n-k}\sum_{j=1}^{n-k}\int_{\ze\in\bC}e^{i\inp{t}{\ze}}\,d\delta_{W_{n,j}(\ome)}(\ze)=\nth{n-k}\sum_{j=1}^{n-k}e^{i\inp{t}{W_{n,j}(\ome)}},\label{eq14}%{eq14}
\end{align}
where $\bP_n(\ome)$ is as in \eqref{eq10} and $\inp{\cdot}{\cdot}$ is the dot product on $\bR^2$.

%%%%% 3.2 PROPOSITION 3 %%%%%
\bigskip
\subsection{Proposition \ref{pro3}}
\label{sec3.2}

\bigskip
We are now ready to prove the `if' part of Proposition \ref{pro3}. Note that the following proof when taking $k=1$ plus \cite[Lemma 2.2]{Sub} and \cite[Proposition 3.2]{Sub} (or Lemma \ref{lem7} below) form an alternative proof of \cite[Theorem 2.3]{Sub}.

%---------- PROOF OF `IF' PART OF PROPOSITION 3 ----------
\bigskip
\begin{proof}[Proof of `if' part of Proposition 4]
{W}e shall verify \eqref{eq13} under assumption \eqref{eq7} and conditions (i) and (ii). From \eqref{eq14} and $\inp{a}{b}=\half\rob{a\ol b+\ol ab}$, we have
\begin{align}
\varphi_{W_n(\ome)}(t)&=\nth{n-k}\sum_{j=1}^{n-k}\sum_{m=0}^\infty\fr{i^m}{2^mm!}\sum_{r=0}^mC^m_rt^{m-r}{\ol t}^r{W_{n,j}(\ome)}^r{\ol{W_{n,j}(\ome)}}^{m-r}\nonumber\\
&=\sum_{m=0}^\infty\fr{i^m}{2^mm!}\sum_{r=0}^mC^m_rt^{m-r}{\ol t}^r\ub{\nth{n-k}\sum_{j=1}^{n-k}{W_{n,j}(\ome)}^r{\ol{W_{n,j}(\ome)}}^{m-r}}_{S_{n,m,r}(\ome)}.\label{eq15}%{eq15}
\end{align}
Write
$$W_{n,j}(\ome):=R_{n,j}(\ome)e^{i\Phi_{n,j}(\ome)},\ R_{n,j}(\ome)\in[0,\infty),\ \Phi_{n,j}(\ome)\in[0,2\pi).$$
By assumption \eqref{eq7}, {\color{b}$\abs{S_{n,m,r}}\leq M^m$} a.s.. This enables us to see that \eqref{eq15} is uniformly absolutely convergent a.s. by {\color{b}checking
$$\sum_{m=0}^\infty\abs{\fr{i^m}{2^mm!}\sum_{r=0}^mC^m_rt^{m-r}{\ol t}^rS_{n,m,r}}\leq\sum_{m=0}^\infty\fr{M^m\abs{t}^m}{2^mm!}\sum_{r=0}^mC^m_r=\sum_{m=0}^\infty\fr{M^m\abs{t}^m}{m!}=e^{M\abs{t}}<\infty.$$
As} a result, if
\begin{equation}\label{eq16}%{eq16}
\lim_{\ntti}S_{n,m,r}
\end{equation}
exists a.s., then we can conclude from \eqref{eq15} that
\begin{equation}\label{eq17}%{eq17}
\varphi_{W_n}(t)\to\sum_{m=0}^\infty\fr{i^m}{2^mm!}\sum_{r=0}^mC^m_rt^{m-r}{\ol t}^r\lim_{\ntti}S_{n,m,r}\antti\text{ \ a.s..}
\end{equation}

Now, we find \eqref{eq16}. For $m-r\geq r$,
\begin{align*}
&S_{n,m,r}(\ome)=\nth{n-k}\sum_{j=1}^{n-k}{R_{n,j}(\ome)}^{2r}{\ol{W_{n,j}(\ome)}}^{m-2r}\\
&\quad=\nth{n-k}\sum_{j=1}^{n-k}{\ol{W_{n,j}(\ome)}}^{m-2r}-\ub{\nth{n-k}\sum_{j=1}^{n-k}\rob{{R_{n,j}(\ome)}^{m-2r}-{R_{n,j}(\ome)}^m}e^{-i(m-2r)\Phi_{n,j}(\ome)}}_{T_{n,m,r}(\ome)}.
\end{align*}
{\color{b}For any arbitrary $0<\rho<1<\rho'$, by assumption \eqref{eq7},
\begin{align*}
\abs{T_{n,m,r}}&\leq\nth{n-k}\sum_{j=1}^{n-k}\abs{{R_{n,j}}^m-{R_{n,j}}^{m-2r}}\\
&=\nth{n-k}\sum_{j=1}^{n-k}\abs{{R_{n,j}}^m-{R_{n,j}}^{m-2r}}\cdot\bIn_{\cub{R_{n,j}\,\leq\,\rho}\,\cup\,\cub{R_{n,j}\,\geq\,\rho'}}(W_{n,j})\\
&\qquad\quad+\nth{n-k}\sum_{j=1}^{n-k}\rob{{R_{n,j}}^{m-2r}-{R_{n,j}}^m}\cdot\bIn_{\cub{\rho\,<\,R_{n,j}\,\leq\,1}}(W_{n,j})\\
&\qquad\quad+\nth{n-k}\sum_{j=1}^{n-k}\rob{{R_{n,j}}^m-{R_{n,j}}^{m-2r}}\cdot\bIn_{\cub{1\,<\,R_{n,j}\,<\,\rho'}}(W_{n,j})\\
&\leq(M^m+M^{m-2r})\cdot\mu_n(\ol{\bD_\rho}\cup({\bD_{\rho'}}^c))\\
&\qquad+(1-\rho^m)\cdot\ub{\mu_n(\ol\bD\setminus\ol{\bD_\rho})}_{\leq1}+\,({\rho'}^m-1)\cdot\ub{\mu_n(\bD_{\rho'}\setminus\ol\bD)}_{\leq1}\text{ \ a.s.}.
\end{align*}
By condition (i), \smash{$\mu_n(\ol{\bD_\rho}\cup({\bD_{\rho'}}^c))\leq\mu_n((\UC)^c)\to0$} as $\ntti$ a.s., so that
$$\limsup_{\ntti}\abs{T_{n,m,r}}\leq{\rho'}^m-\rho^m\text{ \ a.s..}$$
Since both $\rho$ and $\rho'$ can be arbitrarily close to $1$, we indeed have
$$T_{n,m,r}\to0\antti\text{ \ a.s.}.$$
Together} with condition (ii), we get
$$S_{n,m,r}\to\bE\sqb{{\ol Z}^{m-2r}}\antti\text{ \ a.s.}.$$
And for $m-r<r$, a similar argument would give
$$S_{n,m,r}\to\bE\sqb{Z^{2r-m}}=\bE\sqb{{\ol Z}^{m-2r}}\antti\text{ \ a.s..}$$

Consequently, \eqref{eq17} proceeds as
\begin{align}
\varphi_{W_n}(t)&\to\sum_{m=0}^\infty\fr{i^m}{2^mm!}\sum_{r=0}^mC^m_rt^{m-r}{\ol t}^r\bE\sqb{{\ol Z}^{m-2r}}\antti\text{ \ a.s..}\nonumber\\
&=\sum_{m=0}^\infty\bE\sqb{\fr{i^m}{2^mm!}\sum_{r=0}^mC^m_rt^{m-r}{\ol t}^rZ^r{\ol Z}^{m-r}}\nonumber\\
&=\sum_{m=0}^\infty\bE\sqb{\fr{i^m}{m!}\inp{t}{Z}^m}=\bE\sqb{\sum_{m=0}^\infty\fr{i^m}{m!}\inp{t}{Z}^m}=\bE\sqb{e^{i\inp{t}{Z}}}=\varphi_Z(t)\label{eq19}%{eq19}
\end{align}
as desired. (The interchange of infinite summation and expectation in \eqref{eq19} is validated by checking, with the Cauchy--Schwarz inequality, that
$$\sum_{m=0}^\infty\bE\sqb{\abs{\fr{i^m}{m!}\inp{t}{Z}^m}}\leq\sum_{m=0}^\infty\fr{\bE\sqb{\abs{t}^m\abs{Z}^m}}{m!}=\sum_{m=0}^\infty\fr{\abs{t}^m}{m!}=e^{\abs{t}}<\infty,$$
so that a corollary of Lebesgue's dominated convergence theorem (\cite[Theorem 1.38]{Ru}) applies.)
\end{proof}

\bigskip
The above `if' part is sufficient for the rest of this article. We prove the `only if' part for completeness:

%---------- PROOF OF `ONLY IF' PART OF PROPOSITION 3 ----------
\bigskip
\begin{proof}[Proof of `only if' part of Proposition 4]
As remarked before, condition (i) simply follows from Proposition \ref{pro2}, thus it only remains to show condition (ii).

{\color{b}Let $M':=\max\cub{M,1}+\delta$, where $\delta>0$.} For each $p\in\bN$, consider the bounded continuous real-valued {\color{b}function
$$f(\ze):=\left\{\begin{matrix}
\re\,\ze^p&\text{for}&\abs{\ze}\leq M'\vspace{.1cm}\\
\dfr{\re\,\ze^p}{\abs{\ze}^p}&\text{for}&\abs{\ze}>M'
\end{matrix}\right.\text{ \ for \ }\ze\in\bC.$$
Then}, as by assumption \eqref{eq7} {\color{b}and that $M'>M$ so that \smash{$\mu_n(\ol{\bD_{M'}}^c)\leq\mu_n(\ol{\bD_M}^c)=0$}} a.s., using \eqref{eq4} we {\color{b}have
\begin{align}
&\re\rob{\nth{n-k}\sum_{j=1}^{n-k}{W_{n,j}}^p}=\re\rob{\nth{n-k}\sum_{j=1}^{n-k}\int_{\ze\in\bC}\ze^p\,d\delta_{W_{n,j}}(\ze)}=\int_{\ze\in\bC}\re\,\ze^p\,d\mu_n(\ze)\nonumber\\
&\quad=\int_{\ze\in\ol{\bD_{M'}}}\re\,\ze^p\,d\mu_n(\ze)+\int_{\ze\in\ol{\bD_{M'}}^c}\fr{\re\,\ze^p}{\abs{\ze}^p}\,d\mu_n(\ze)=\int_{\ze\in\bC}f(\ze)\,d\mu_n(\ze)\text{ \ a.s..}\label{eq20}%{eq20}
\end{align}
On} the other hand, as $\mu$ is supported in $\UC$ {\color{b}and that $M'>1$} so that {\color{b}\smash{$\mu(\ol{\bD_{M'}}^c)\leq$} \smash{$\mu((\UC)^c)=0$}}, using \eqref{eq1} we {\color{b}have
\begin{align}
\re\,\bE\sqb{Z^p}&=\bE\sqb{\re\,Z^p}=\int_{\ze\in\bC}\re\,\ze^p\,d\mu(\ze)\nonumber\\
&=\int_{\ze\in\ol{\bD_{M'}}}\re\,\ze^p\,d\mu(\ze)+\int_{\ze\in\ol{\bD_{M'}}^c}\fr{\re\,\ze^p}{\abs{\ze}^p}\,d\mu(\ze)=\int_{\ze\in\bC}f(\ze)\,d\mu(\ze).\label{eq21}%{eq21}
\end{align}
Therefore}, by \eqref{eq20} and \eqref{eq21}, $\mu_n\wto\mu$ as $\ntti$ a.s. implies that
$$\begin{matrix}
\dst\re\rob{\nth{n-k}\sum_{j=1}^{n-k}{W_{n,j}}^p}=\int_{\ze\in\bC}f(\ze)\,d\mu_n(\ze)\to\int_{\ze\in\bC}f(\ze)\,d\mu(\ze)=\re\,\bE\sqb{Z^p}\\
\antti\text{ \ a.s..}
\end{matrix}$$
A similar argument with `$\re$' in the construction of $f$ replaced by `$\im$' would give us
\begin{equation*}
\im\rob{\nth{n-k}\sum_{j=1}^{n-k}{W_{n,j}}^p}\to\im\,\bE\sqb{Z^p}\antti\text{ \ a.s.}.\qedhere
\end{equation*}
\end{proof}

%%%%% 3.3 A KEY LEMMA %%%%%
\bigskip
\subsection{A Key Lemma}
\label{sec3.3}

\bigskip
The lemma central to the rest of this article is contained in \cite[Proposition 3.2]{Sub}, in which the proof uses a result for polynomials in \cite{CheuNg1} and \cite{KoRi} obtained by a companion matrix approach. We present {\color{b}essentially the} same proof, but somewhat more straightforward:

%---------- LEMMA 7 ----------
\bigskip
\begin{lem}\label{lem7}
Let $p\in\bN$. Then, the mean
$$\fr{{w_1}^p+\cdots+{w_{n-1}}^p}{n-1}$$
of the {\color{b}$p$-th powers} of the zeros $w_1,\dots,w_{n-1}$ of $Q$ in \eqref{eq9} can be expressed (in terms of those zeros $z_j$ of $P$ in \eqref{eq8} and the constants $\lam_j$) as
\begin{align}
\fr{n}{n-1}\fr{{z_1}^p+\cdots+{z_n}^p}{n}&-\fr{p}{n-1}\sum_{j=1}^n\al_j{z_j}^p\nonumber\\
&+\nth{n-1}\,{\sum}'(-1)^s\prod_{t=1}^{s-1}{\color{b}\rob{\sum_{j=1}^n\al_j{z_j}^{h_t}}}\sum_{j=1}^n\al_j{z_j}^{q+r},\label{eq22}%{eq22}
\end{align}
where
\begin{equation}\label{eq23}%{eq23}
\al_j:=\fr{\lam_j}{\dst\sum_{j=1}^n\lam_j}\text{ \ and \ }{\sum}':=\sum_{q=1}^{p-1}\sum_{r=0}^{p-q-1}\sum_{s=2}^{p-q-r+1}\sum_{\substack{h_1,\dots,h_{s-1}\in\bN\\h_1+\dots+h_{s-1}=p-q-r}}.
\end{equation}
\end{lem}

%---------- PROOF OF LEMMA 7 ----------
\bigskip
\begin{proof}
We shall use \cite[Theorem 1.2]{CheuNg2} (instead of the less general results in \cite{CheuNg1} and \cite{KoRi}), which asserts that:
\begin{quote}
\it If $D$ is an $n\times n$ matrix with characteristic polynomial $\wt P(z):=(z-z_1)\cdots(z-z_n)$ and $\wt Q(z)$ is a monic polynomial of degree $n-1$ given by
$$\fr{\wt Q(z)}{\wt P(z)}:=\sum_{j=1}^n\fr{\al_j}{z-z_j}\text{ \ (as in \eqref{eq9})},$$
then there exists a rank one matrix $H$ such that $H^2=H$ and $z\wt Q(z)$ is the characteristic polynomial of $D-DH$. In particular, if $D=\diag(z_1,\dots,z_n)$, then we can take $H=LJ$, where $L=\diag(\al_1,\dots,\al_n)$ and $J$ is the $n\times n$ matrix all of whose entries are one.
\end{quote}
Taking the two monic polynomials in this cited theorem as
$$\wt P=P\quad\text{and}\quad\wt Q=\fr{Q(z)}{\dst\sum_{j=1}^n\lam_j}=P(z)\sum_{j=1}^n\fr{\al_j}{z-z_j},$$
where $P$ and $Q$ are as in \eqref{eq8} and \eqref{eq9} and $\al_j$ is as in \eqref{eq23}, we know that $0$, $w_1,\dots,w_{n-1}$ are the eigenvalues of
$$D-DLJ,$$
where $D$, $L$ and $J$ are as just mentioned. Since
$$0^p+{w_1}^p+\cdots+{w_{n-1}}^p=\tr\rob{(D-DLJ)^p},$$
we expand the matrix power $(D-DLJ)^p$ and obtain
$$D^p+\sum_{q=1}^pD^q(-LJ)D^{p-q}+{\sum}'D^q(-LJ)\ub{D^{h_1}(-LJ)\cdots D^{h_{s-1}}(-LJ)}_{s-1\text{ factors of the form }D^{h_\cdot}(-LJ)}D^r,$$
where ${\sum}'$ is as in \eqref{eq23}. This expansion can be verified by induction  on $p$. Then, using the special property
$$JD^hLJ=\tr(D^hL)J\ \ \forall\ h\in\bN$$
(as is easy to check) repeatedly, the expansion becomes
\begin{align}
&\quad\ D^p-\sum_{q=1}^pD^qLJD^{p-q}+{\sum}'(-1)^sD^qL(JD^{h_1}LJ)\ub{D^{h_2}(LJ)\cdots D^{h_{s-1}}(LJ)}_{s-2\text{ factors remain}}D^r\nonumber\\
&=D^p-\sum_{q=1}^pD^qLJD^{p-q}\nonumber\\
&\qquad+{\sum}'(-1)^s\tr(D^{h_1}L)D^qL(JD^{h_2}LJ)\ub{D^{h_3}(LJ)\cdots D^{h_{s-1}}(LJ)}_{s-3\text{ factors remain}}D^r\nonumber\\
&=\cdots=D^p-\sum_{q=1}^pD^qLJD^{p-q}+{\sum}'(-1)^s\tr(D^{h_1}L)\cdots\tr(D^{h_{s-1}}L)D^qLJD^r.\label{eq24}%{eq24}
\end{align}
Finally, the facts
$$\tr(D^hL)=\sum_{j=1}^n\al_j{z_j}^h\text{ \ and \ }\tr(D^hLJD^{h'})=\sum_{j=1}^n\al_j{z_j}^{h+h'}\ \ \forall\ h,h'\in\bN\cup\cub{0}$$
(as are also easy to check) show that the trace of \eqref{eq24} actually equals
\begin{equation*}{z_1}^p+\cdots+{z_n}^p-\sum_{q=1}^p\sum_{j=1}^n\al_j{z_j}^p+{\sum}'(-1)^s\sum_{j=1}^n\al_j{z_j}^{h_1}\cdots\sum_{j=1}^n\al_j{z_j}^{h_{s-1}}\sum_{j=1}^n\al_j{z_j}^{q+r}.\qedhere
\end{equation*}
\end{proof}

%%%%% 3.4 COROLLARY 4%%%%%
\bigskip
\subsection{Corollary \ref{cor4}}
\label{sec3.4}

\bigskip
Now, $W_{n,1}(\ome),\dots,W_{n,n-k}(\ome)$ in \eqref{eq3} are the zeros of the {\color{b}$k$-th} order derivative of $P_n(\ome)$ in \eqref{eq2}, thus condition \eqref{eq7} directly follows from {\color{b}a repeated use of} the Gauss--Lucas theorem, and it only remains to verify conditions (i) and (ii) in Proposition \ref{pro3}, so that \smash{\color{b}$\mu^{(k)}_n\wto\mu$} as $\ntti$ a.s.:

%---------- FIRST PART OF PROOF OF COROLLARY 4 ----------
\bigskip
\begin{proof}[Verification of condition (i) in Proposition 4]
Since all \smash{\color{b}$|W^{(k)}_{n,j}|\leq1$} a.s. so that \smash{\color{b}$\mu^{(k)}_n$} is supported in $\ol\bD$ a.s., it suffices to show that for any $0<r<1$,
$${\color{b}\mu^{(k)}_n}(\ol{\bD_r})\to0\antti\text{ \ a.s.,}$$
where $\bD_r:=\cub{z\in\bC:\abs{z}<r}$ denotes the disc of radius $r$.

Let $0<r<1$. Since all $\abs{Z_j}=1$ a.s. and $\mu$ is non-uniform, let as in the proof of \cite[Lemma 2.2]{Sub}
\begin{align*}
U_{n,1}(z)&:=\nth{n}\fr{{P_n}'(z)}{P_n(z)}=\nth{n}\sum_{j=1}^n\nth{z-Z_j}\\
&\phantom{:}=\sum_{m=0}^\infty\nth{n}\sum_{j=1}^n-\ol{Z_j}^{m+1}z^m\to F(z):=\sum_{m=0}^\infty-\ol{a_m}z^m\not\equiv0\antti
\end{align*}
uniformly on compact subsets of $\bD$ a.s., where
$$a_m:=\bE\sqb{Z^{m+1}}\text{ \ for \ }m=0,1,\dots.$$
To deal with higher order derivatives of $P_n$, we employ the trick
$$\nth{n^l}\fr{P^{(l)}(z)}{P(z)}=\nth{n}\rob{\nth{n^{l-1}}\fr{P^{(l-1)}(z)}{P(z)}}'+\nth{n}\fr{P'(z)}{P(z)}\cdot\nth{n^{l-1}}\fr{P^{(l-1)}(z)}{P(z)}\ \ \forall\ l\in\bN.$$
Indeed, with the facts about $U_{n,1}$ and $F$ above, we have
\begin{align*}
U_{n,2}(z)&:=\nth{n^2}\fr{{P_n}''(z)}{P_n(z)}=\nth{n}{U_{n,1}}'(z)+U_{n,1}(z)\cdot U_{n,1}(z)\to0+F(z)^2,\\
U_{n,3}(z)&:=\nth{n^3}\fr{{P_n}^{(3)}(z)}{P_n(z)}=\nth{n}{U_{n,2}}'(z)+U_{n,1}(z)\cdot U_{n,2}(z)\to0+F(z)^3,\\
&\qquad\qquad\qquad\qquad\qquad\qquad\vdots\\
U_{n,k}(z)&:=\nth{n^k}\fr{{P_n}^{(k)}(z)}{P_n(z)}=\nth{n}{U_{n,k-1}}'(z)+U_{n,1}(z)\cdot U_{n,k-1}(z)\to0+F(z)^k\\
&\qquad\qquad\qquad\qquad\quad\antti
\end{align*}
uniformly on compact subsets of $\bD$ a.s.. Now, $U_{n,k}$ is also a holomorphic function on $\bD$ whose zeros are exactly those of ${P_n}^{(k)}$ a.s.. Pick any $r<r'<1$ such that $F$ has no zero on $\partial\bD_{r'}$. Since $F^k\not\equiv0$ is also holomorphic on $\bD$, it has only finitely many, say $M$, zeros in $\ol{\bD_{r'}}$. And as $F^k$ has no zero on $\partial\bD_{r'}$, it then follows from Hurwitz's theorem (\cite[Theorem 2.5]{Co}) that $U_{n,k}$, as well as ${P_n}^{(k)}$, also has exactly $M$ zeros in $\ol{\bD_{r'}}$ for sufficiently large $n$ a.s.. As a result,
\begin{equation*}
{\color{b}\mu^{(k)}_n}(\ol{\bD_r})=\nth{n-k}\sum_{j=1}^{n-k}\bIn_{\ol{\bD_r}}(W_{n,j})\leq\fr{M}{n-k}\to0\antti\text{ \ a.s.}.\qedhere
\end{equation*}
\end{proof}

%---------- SECOND PART OF PROOF OF COROLLARY 4 ----------
\bigskip
\begin{proof}[Verification of condition (ii) in Proposition 4]
\ Recall that when $\lam_j=1$ in \eqref{eq9}, we would have $Q=P'$ and $\al_j=\nth{n}$ in Lemma \ref{lem7}. Also, if $\abs{z_j}=1$, we have
$$\abs{\sum_{j=1}^n\al_j{z_j}^h}\leq\sum_{j=1}^n\abs{\al_j}=1\ \ \forall\ h\in\bN.$$
As a result, the second term and thereafter in \eqref{eq22} is merely $O\rob{\nth{n}}$ (which may depend on $p$), and then we can write the result of Lemma \ref{lem7} as
$$\fr{{w_1}^p+\cdots+{w_{n-1}}^p}{n-1}=\fr{{z_1}^p+\cdots+{z_n}^p}{n}+O\rob{\nth{n}}.$$
Applying Lemma \ref{lem7} in this form $k$ times respectively with
$$P={P_n},{P_n}',\dots,{P_n}^{(k-1)}\text{ \ and \ }Q={P_n}',{P_n}'',\dots,{P_n}^{(k)},$$
we have
\begin{align*}
\fr{{W^{(k)}_{n,1}}^p+\cdots+{W^{(k)}_{n,n-k}}^p}{n-k}&=\fr{{W^{(k-1)}_{n,1}}^p+\cdots+{W^{(k-1)}_{n,n-k+1}}^p}{n-k+1}+O\rob{\nth{n}}\\
&\qquad\qquad\qquad\quad\vdots\\
&=\fr{{W^{(1)}_{n,1}}^p+\cdots+{W^{(1)}_{n,n-1}}^p}{n-1}+O\rob{\nth{n}}\\
&=\fr{{Z_1}^p+\cdots+{Z_n}^p}{n}+O\rob{\nth{n}}\to\bE\sqb{Z^p}\antti\text{ \ a.s.}
\end{align*}
by Kolmogorov's strong law of large numbers, where \smash{$W^{(l)}_{n,1},\dots,W^{(l)}_{n,n-l}$} denote the zeros of the {\color{b}$l$-th} order derivative of $P_n$.
\end{proof}

%%%%% 3.5 COROLLARY 5 %%%%%
\bigskip
\subsection{Corollary \ref{cor5}}
\label{sec3.5}

\bigskip
Under the assumptions in Corollary \ref{cor5}, we verify conditions (i) and (ii) in Proposition \ref{pro3}, so that $\mu_n\wto\mu$ as $\ntti$ a.s.:

%---------- FIRST PART OF PROOF OF COROLLARY 5 ----------
\bigskip
\begin{proof}[Verification of condition (i) in Proposition 4]
\ Since, by assumption, all $\abs{W_{n,j}}\leq1$ a.s. so that $\mu_n$ is supported in $\ol\bD$, it suffices to show that for any $0<r<1$,
$$\mu_n(\ol{\bD_r})\to0\antti\text{ \ a.s..}$$

Let $0<r<1$. Since all $\abs{Z_j}=1$ a.s., for $z\in\bD$,
\begin{align}
V_n(z)&:=\nth{n}\fr{Q_n(z)}{P_n(z)}=\nth{n}\sum_{j=1}^n\fr{\lam_{n,j}}{z-Z_j}=\nth{n}\sum_{j=1}^n\fr{-\lam_{n,j}\ol{Z_j}}{1-\ol{Z_j}z}\nonumber\\
&\phantom{:}=\nth{n}\sum_{j=1}^n-\lam_{n,j}\ol{Z_j}\sum_{m=0}^\infty\ol{Z_j}^mz^m=\sum_{m=0}^\infty\ub{\nth{n}\sum_{j=1}^n-\lam_{n,j}\ol{Z_j}^{m+1}}_{B_{n,m}}z^m\label{eq25}%{eq25}
\end{align}
is a holomorphic function on $\bD$ a.s. whose zeros are exactly those of $Q_n(z)$. By conditions (i) and (ii), we note that $\abs{b_m}\leq K+1$ for sufficiently large $m$, so for $z\in\bD$,
$$G(z):=\sum_{m=0}^\infty-\ol{b_m}z^m\not\equiv0$$
is also a holomorphic function on $\bD$ ($G\not\equiv0$ comes from the assumption that at least one $b_m\neq0$). Pick any $r<r'<1$ such that $G$ has no zero on $\partial\bD_{r'}$. It is clear from conditions (i) and (ii) that
$$V_n(z)=\sum_{m=0}^\infty B_{n,m}z^m$$
is uniformly absolutely convergent on $\ol{\bD_{r'}}$ a.s., so
\begin{equation}\label{eq26}%{eq26}
V_n(z)\to\sum_{m=0}^\infty\lim_{\ntti}B_{n,m}z^m=\sum_{m=0}^\infty-\ol{b_m}z^m=G(z)\antti
\end{equation}
uniformly on $\ol{\bD_{r'}}$ a.s.. And as $G$ has no zero on $\partial\bD_{r'}$, it then follows from Hurwitz's theorem that $V_n$, as well as $Q_n$, also has exactly the same finite number, say $M$, of zeros in $\ol{\bD_{r'}}$ for sufficiently large $n$ a.s.. As a result,
\begin{equation*}
\mu_n(\ol{\bD_r})=\nth{n-k}\sum_{j=1}^{n-k}\bIn_{\ol{\bD_r}}(W_{n,j})\leq\fr{M}{n-k}\to0\antti\text{ \ a.s.}.\qedhere
\end{equation*}
\end{proof}

%---------- SECOND PART OF PROOF OF COROLLARY 5 ----------
\bigskip
\begin{proof}[Verification of condition (ii) in Proposition 4]
For each $p\in\bN$, applying Lemma \ref{lem7}, we have
\begin{align}
\fr{{W_{n,1}}^p+\cdots+{W_{n,n-1}}^p}{n-1}&=\fr{n}{n-1}\fr{{Z_1}^p+\cdots+{Z_n}^p}{n}-\fr{p}{n-1}\sum_{j=1}^n\al_{n,j}{Z_j}^p\nonumber\\
&\qquad+\nth{n-1}\,{\sum}'(-1)^s\prod_{t=1}^{s-1}{\color{b}\rob{\sum_{j=1}^n\al_{n,j}{Z_j}^{h_t}}}\sum_{j=1}^n\al_{n,j}{Z_j}^{q+r},\label{eq27}%{eq27}
\end{align}
a.s., where ${\sum}'$ is as in \eqref{eq23} and
$$\al_{n,j}:=\fr{\lam_{n,j}}{\dst\sum_{j=1}^n\lam_{n,j}}.$$
For each $h\in\bN$, by conditions (ii) and (iii) we have
$$\abs{\sum_{j=1}^n\al_{n,j}{Z_j}^h}\leq\fr{\ \dst\nth{n}\sum_{j=1}^n\abs{\lam_{n,j}}\ }{\dst\nth{n}\abs{\sum_{j=1}^n\lam_{n,j}}}<\fr{K}{\delta}\ \ \forall\ n\text{ \ a.s..}$$
As a result, the second term and thereafter in \eqref{eq27} is merely $O\rob{\nth{n}}$ (which may depend on $p$) a.s., and then
$$\fr{{W_{n,1}}^p+\cdots+{W_{n,n-1}}^p}{n-1}=\fr{{Z_1}^p+\cdots+{Z_n}^p}{n}+O\rob{\nth{n}}\to\bE\sqb{Z^p}\antti\text{ \ a.s.}$$
by Kolmogorov's strong law of large numbers.
\end{proof}

%---------- COROLLARY 5 -- POLAR DERIVATIVE ----------
\bigskip
\begin{proof}[Polar derivative]
In the case that $Q_n$ is the polar derivative of $P_n$ with respect to {\color{b}any $\xi\in\bC$ with $\abs{\xi}>1$}, conditions (ii) and (iii)
$$\nth{n}\sum_{j=1}^n\abs{\lam_{n,j}}\leq\abs{\xi}+1\text{ \ and \ }\nth{n}\abs{\sum_{j=1}^n\lam_{n,j}}\geq\abs{\xi}-1>0$$
hold because all $\abs{Z_j}=1$ a.s.. Moreover, since $\abs{\xi}>1$, condition \eqref{eq7} then follows from Laguerre's theorem (\cite[Theorem A]{Az}, \cite[Lemma 1.1.11]{Pr}, \cite[Theorem 3.2.1(i), Theorem 3.2.1a]{RahSc}, which can be regarded as the Gauss--Lucas theorem for polar derivative). For condition (i), we actually have
\begin{align}
b_m=\lim_{\ntti}\nth{n}\sum_{j=1}^n\ol{\lam_{n,j}}{Z_j}^{m+1}&=\ol\xi\lim_{\ntti}\nth{n}\sum_{j=1}^n{Z_j}^{m+1}-\lim_{\ntti}\nth{n}\sum_{j=1}^n{Z_j}^m\text{ \ a.s.}\nonumber\\
&=\ol\xi\,\bE\sqb{Z^{m+1}}-\bE\sqb{Z^m}\label{eq28}%{eq28}
\end{align}
by Kolmogorov's strong law of large numbers. In {\color{b}particular,
$$b_0=\ol\xi\,\bE\sqb{Z}-1\left\{\begin{alignedat}{4}
&=-1									&&\text{ \ if \ }\mu_{\UC}\text{ \ is uniform,}\\
&\neq0\text{ \ for \ }\abs{\xi}>\nth{\abs{\bE[Z]}}	&&\text{ \ if \ }\mu_{\UC}\text{ \ is non-uniform.}
\end{alignedat}\right.$$
Here, we have used \cite[Lemma 3.1]{Sub} that `$\mu_{\UC}$ is uniform' if and only if $\bE\sqb{Z^p}=0$ for all $p\in\bN$.}
\end{proof}

%---------- REMARK 1 ----------
\bigskip
\begin{remark}
This result for polar derivative is independent of whether $\mu$ is uniform or not as long as {\color{b}$\abs{\xi}\gg1$}. Such flexibility and \eqref{eq25} are the key of the proof of Corollary \ref{cor6}.
\end{remark}

%---------- COROLLARY 5 -- SZ.-NAGY'S GENERALIZED DERIVATIVE ----------
\bigskip
\begin{proof}[Sz.-Nagy's generalized derivative]
In the case that $Q_n$ is a Sz.-Nagy's generalized derivative of $P_n$, we have
$$\nth{n}\sum_{j=1}^n\abs{\lam_{n,j}}=\nth{n}\sum_{j=1}^n\lam_{n,j}=1,$$
so conditions (ii) and (iii) hold. Furthermore, by the Gauss--Lucas theorem for Sz.-Nagy's generalized derivative (\cite[Corollary 4.1]{CheuNg2} (proved by matrix methods); also mentioned in \cite[p.115]{RahSc}), we obtain condition \eqref{eq7}. Finally, for condition (i), since $Z=Z_1,Z_2,\dots$ are i.i.d. and all {\color{b}the} degenerate {\color{b}random variables} $\lam_{n,j}\leq M$ a.s., we apply \cite[Theorem 5]{ChoSun1} with
$$X_j={Z_j}^{m+1}-\bE\sqb{Z^{m+1}}\text{ \ and \ }a_{n,j}=\fr{\lam_{n,j}}{n}$$
to get
$$\nth{n}\sum_{j=1}^n\lam_{n,j}{Z_j}^{m+1}-\bE\sqb{Z^{m+1}}=\sum_{j=1}^na_{n,j}X_j\to0\antti\text{ \ a.s.};$$
but $\mu$ is non-uniform, so at least one of $b_m:=\bE\sqb{Z^{m+1}}\neq0$.
\end{proof}

%---------- REMARK 2 ----------
\bigskip
\begin{remark}
For the case of Sz.-Nagy's generalized derivative, besides the degeneracy and a.s. uniform boundedness of $\lam_{n,j}$, there are plenty of sufficient conditions for condition (i) to hold. See, for instance, \cite{ChenGa}, \cite{ChoSun2}, \cite{Cu} and \cite{Sun}.
\end{remark}

%%%%% 3.6 COROLLARY 6 %%%%%
\bigskip
\subsection{Corollary \ref{cor6}}
\label{sec3.6}

\bigskip
For each $l\in\bN$, let
\begin{equation}\label{eq29}%{eq29}
W^{(l)}_{n,1}(\ome),\dots,W^{(l)}_{n,n-l}(\ome)
\end{equation}
(as in \eqref{eq3}) be the zeros of the {\color{b}$l$-th} order derivative of $P_n(\ome)$ in \eqref{eq2}. By the fact that $\mu$ is supported in $\UC$ and the Gauss--Lucas theorem, \smash{$|W^{(l)}_{n,j}|\leq1$} for all $n$ and $j$ a.s.. From the second part of the proof of Corollary \ref{cor4}, we see that condition (ii) in Proposition \ref{pro3} holds automatically for \eqref{eq29}. Thus in this case, Proposition \ref{pro3} reduces to simply
\begin{equation}\label{eq30}%{eq30}
\mu^{(l)}_n\wto\mu\antti\text{ \ a.s.}\quad\iff\quad\mu^{(l)}_n\arr\UC\antti\text{ \ a.s.},
\end{equation}
where $\mu^{(l)}_n$ denotes the empirical measure of \eqref{eq29}. We shall prove Corollary \ref{cor6} by establishing:
$$\text{\it For each \ }0<r<1,\ \mu^{(l)}_n(\ol{\bD_r})=0\text{ \ \it for sufficiently large \ }n\text{ \ \it a.s..}$$
inductively on $l$. Now, since $\mu$ is supported in $\UC$ and uniform, by \cite[Lemma 3.1]{Sub} (as in the first part of the proof of Corollary \ref{cor4}),
\begin{equation}\label{eq31}%{eq31}
\bE\sqb{Z^p}=0\ \ \forall\ p\in\bN.
\end{equation}

%---------- COROLLARY 6 -- FIRST ORDER DERIVATIVE ----------
\bigskip
\begin{proof}[First order derivative]
Let $0<r<1$. Taking \eqref{eq28} and \eqref{eq31} into account, \eqref{eq25} and \eqref{eq26} with $\lam_{n,j}=\xi-Z_j$ now say that
$$\nth{n}\fr{D_\xi P_n(z)}{P_n(z)}\to1\antti$$
uniformly on $\ol{\bD_r}$ for any {\color{b}$\xi\in\bC$} a.s.. Then by Hurwitz's theorem, $D_\xi P_n$ has no zero in $\ol{\bD_r}$, i.e.
$$\mu^\xi_n(\ol{\bD_r})=0,$$
for sufficiently large $n$ a.s., where $\mu^\xi_n$ denotes the empirical measure of the zeros of $D_\xi P_n$. For these $n$, since
$$\fr{D_\xi P_n}{\xi}\to{P_n}'\text{ \ as \ }\xi\to\infty$$
uniformly on $\ol{\bD_r}$ (\cite[p.185]{Shei}) a.s., we further deduce by a corollary of Hurwitz's theorem (\cite[Corollary 2.6]{Co}) that ${P_n}'$ also has no zero in $\ol{\bD_r}$, i.e.
\begin{equation*}
\mu^{(1)}_n(\ol{\bD_r})=0\text{ \ a.s..}\qedhere
\end{equation*}
\end{proof}

%---------- COROLLARY 6 -- INDUCTION ----------
\bigskip
\begin{proof}[Induction]
Let $l\geq2$. Assume that for each $0<r<1$,
\begin{equation}\label{eq32}%{eq32}
\mu^{(l-1)}_n(\ol{\bD_r})=0\text{ \ for sufficiently large \ }n\text{ \ a.s.},
\end{equation}
so that ($\mu^{(l-1)}_n(\ol{\bD_r})\to0$ as $\ntti$ a.s., and then)
\begin{equation}\label{eq33}%{eq33}
\mu^{(l-1)}_n\wto\mu\antti\text{ \ a.s.}.
\end{equation}

Let $r_0<r_0'<1$. By induction hypothesis \eqref{eq32}, for sufficiently large $n$,
\begin{equation}\label{eq34}%{eq34}
\abs{W^{(l-1)}_{n,j}}>r_0'\ \ \forall\ j\text{ \ a.s.},
\end{equation}
so that for $z\in\ol{\bD_{r_0'}}$,
\begin{align}
&H^{(l-1),\xi}_n(z):=\nth{n-l+1}\fr{D_\xi{P_n}^{(l-1)}(z)}{{P_n}^{(l-1)}(z)}=\nth{n-l+1}\sum_{j=1}^{n-l+1}\fr{\xi-W^{(l-1)}_{n,j}}{z-W^{(l-1)}_{n,j}}\nonumber\\
&\qquad=\nth{n-l+1}\sum_{j=1}^{n-l+1}\fr{1-\dfr{\xi}{W^{(l-1)}_{n,j}}}{\raisebox{-.25cm}{$1-\dfr{z}{W^{(l-1)}_{n,j}}$}}=\nth{n-l+1}\sum_{j=1}^{n-l+1}\rob{1-\fr{\xi}{W^{(l-1)}_{n,j}}}\sum_{m=0}^\infty\fr{z^m}{{W^{(l-1)}_{n,j}}^m}\nonumber\\
&\qquad=\sum_{m=0}^\infty\ub{\nth{n-l+1}\sum_{j=1}^{n-l+1}\rob{\nth{{W^{(l-1)}_{n,j}}^m}-\fr{\xi}{{W^{(l-1)}_{n,j}}^{m+1}}}}_{C^{(l-1),\xi}_{n,m}}z^m\label{eq35}%{eq35}
\end{align}
is a holomorphic function on $\ol{\bD_{r_0'}}$ a.s. for any {\color{b}$\xi\in\bC$}. We shall find
$$\lim_{\ntti}C^{(l-1),\xi}_{n,m}.$$

For each $p\in\bN$, consider the bounded continuous real-valued function
$$g(\ze):=\left\{\begin{matrix}
0&\text{for}&\ze=0\vspace{.1cm}\\
\dst\fr{\abs{\ze}^{p+1}}{{r_0'}^{p+1}}\re\,\nth{\ze^p}&\text{for}&0<\abs{\ze}\leq r_0'\vspace{.1cm}\\
\dst\re\,\nth{\ze^p}&\text{for}&\abs{\ze}>r_0'
\end{matrix}\right.\text{ \ for \ }\ze\in\bC.$$
By \eqref{eq34} so that $\mu^{(l-1)}_n$ is supported in $\cub{z\in\bC:r_0'<\abs{z}\leq1}$ a.s., we have
\begin{align}
&\quad\re\rob{\nth{n-l+1}\sum_{j=1}^{n-l+1}\nth{{W^{(l-1)}_{n,j}}^p}}\nonumber\\
&=\re\rob{\nth{n-l+1}\sum_{j=1}^{n-l+1}\int_{\ze\in\bC}\nth{\ze^p}\,d\delta_{W^{(l-1)}_{n,j}}(\ze)}=\int_{\ze\in\bC}\re\,\nth{\ze^p}\,d\mu_n^{(l-1)}(\ze)\nonumber\\
&=\int_{\ze\in\ol{\bD_{r_0'}}\setminus\cub{0}}\fr{\abs{\ze}^{p+1}}{{r_0'}^{p+1}}\re\,\nth{\ze^p}\,d\mu_n^{(l-1)}(\ze)+\int_{\ze\in\ol\bD\setminus\ol{\bD_{r_0'}}}\re\,\nth{\ze^p}\,d\mu_n^{(l-1)}(\ze)\nonumber\\
&\qquad\qquad\qquad\quad+\int_{\ze\in\ol\bD^c}\re\,\nth{\ze^p}\,d\mu_n^{(l-1)}(\ze)=\int_{\ze\in\bC}g(\ze)\,d\mu_n^{(l-1)}(\ze)\text{ \ a.s.}.\label{eq36}%{eq36}
\end{align}
On the other hand, as $\mu$ is supported in $\UC$ so that $\mu((\UC)^c)=0$, using \eqref{eq1} we have
\begin{align}
\re\,\bE\sqb{\nth{Z^p}}&=\bE\sqb{\re\,\nth{Z^p}}=\int_{\ze\in\bC}\re\,\nth{\ze^p}\,d\mu(\ze)\nonumber\\
&=\int_{\ze\in\ol{\bD_{r_0'}}\setminus\cub{0}}\fr{\abs{\ze}^{p+1}}{{r_0'}^{p+1}}\re\,\nth{\ze^p}\,d\mu(\ze)+\int_{\ze\in\bD\setminus\ol{\bD_{r_0'}}}\re\,\nth{\ze^p}\,d\mu(\ze)\nonumber\\
&\qquad+\int_{\ze\in\UC}\re\,\nth{\ze^p}\,d\mu(\ze)+\int_{\ze\in\ol\bD^c}\re\,\nth{\ze^p}\,d\mu(\ze)=\int_{\ze\in\bC}g(\ze)\,d\mu(\ze).\label{eq37}%{eq37}
\end{align}
Therefore, by \eqref{eq36} and \eqref{eq37}, \eqref{eq33} implies that
\begin{align*}
\re\rob{\nth{n-l+1}\sum_{j=1}^{n-l+1}\nth{{W^{(l-1)}_{n,j}}^p}}&=\int_{\ze\in\bC}g(\ze)\,d\mu_n^{(l-1)}(\ze)\\
&\to\int_{\ze\in\bC}g(\ze)\,d\mu(\ze)=\re\,\bE\sqb{\nth{Z^p}}\antti\text{ \ a.s.}.
\end{align*}
A similar argument with `$\re$' in the construction of $g$ replaced by `$\im$' would give us
$$\im\rob{\nth{n-l+1}\sum_{j=1}^{n-l+1}\nth{{W^{(l-1)}_{n,j}}^p}}\to\im\,\bE\sqb{\nth{Z^p}}\antti\text{ \ a.s.}.$$
The overall result is, by \eqref{eq31},
$$\nth{n-l+1}\sum_{j=1}^{n-l+1}\nth{{W^{(l-1)}_{n,j}}^p}\to\bE\sqb{\nth{Z^p}}=\bE\sqb{{\ol Z}^p}=0\antti\text{ \ a.s.,}$$
thus in \eqref{eq35}
$$C^{(l-1),\xi}_{n,m}\to\left\{\begin{matrix}
1\antti&\text{for}&m=0\\
0\antti&\text{for}&m\geq1
\end{matrix}\right.\text{ \ a.s.,}$$
so that by the obvious uniform absolute convergence of \eqref{eq35}, we have
$$H^{(l-1),\xi}_n(z)\to\sum_{m=0}^\infty\lim_{\ntti}C^{(l-1),\xi}_{n,m}z^m=1\antti$$
uniformly on $\ol{\bD_{r_0'}}$ for any {\color{b}$\xi\in\bC$} a.s.. Then by Hurwitz's theorem, $D_\xi{P_n}^{(l-1)}$ has no zero in $\ol{\bD_{r_0}}$, i.e.
$$\mu^{(l-1),\xi}_n(\ol{\bD_{r_0}})=0,$$
for sufficiently large $n$ a.s., where $\mu^{(l-1),\xi}_n$ denotes the empirical measure of the zeros of $D_\xi{P_n}^{(l-1)}$. For these $n$, since
$$\fr{D_\xi{P_n}^{(l-1)}}{\xi}\to({P_n}^{(l-1)})'={P_n}^{(l)}\text{ \ as \ }\xi\to\infty$$
uniformly on $\ol{\bD_{r_0'}}$, we further deduce, by the corollary of Hurwitz's theorem just mentioned, that ${P_n}^{(l)}$ also has no zero in $\ol{\bD_{r_0}}$, i.e.
$$\mu^{(l)}_n(\ol{\bD_{r_0}})=0\text{ \ a.s..}$$
Hence, by \eqref{eq30} and induction, $\mu^{(k)}_n\wto\mu$ as $\ntti$ a.s..
\end{proof}

%%%%%%%%%% 4. CLOSING REMARK %%%%%%%%%%
\bigskip
\section{Closing remark}
\label{sec4}

\bigskip
Recall that the metric topology of weak convergence on the space $\bM$ of probability measures on $(\bC,\cB)$ (mentioned in Section \ref{sec1}) is given by the Prohorov metric $\pi$ (\cite[Theorem 6.8]{Bi}): For $m',m''\in\bM$,
\begin{equation}\label{eq38}%{eq38}
\begin{array}{rl}
\pi(m',m''):=\inf\{\vep>0:&m'(A)\leq m''(A^\vep)+\vep\text{ \ and}\\
&m''(A)\leq m'(A^\vep)+\vep\text{ \ for all \ }A\in\cB\},
\end{array}
\end{equation}
where
$$A^\vep:=\bigcup_{a\in A}\bD_\vep(a)\text{ \ and \ }\bD_\vep(a):=\cub{z\in\bC:\abs{z-a}<\vep}.$$
Consider the case of empirical measures:
\begin{equation}\label{eq39}%{eq39}
m'=\nth{n'}\sum_{i=1}^{n'}\delta_{a'_i}\quad\text{and}\quad m''=\nth{n''}\sum_{j=1}^{n''}\delta_{a''_j}.
\end{equation}
Let $0<q<1$ be fixed. Assume that $m'$ and $m''$ are close with respect to the Prohorov metric $\pi$, say,
\begin{equation}\label{eq40}%{eq40}
\pi(m',m'')<\vep_0\text{ \ for some \ }0<\vep_0\leq1-q.
\end{equation}
Observe that by \eqref{eq38} and assumption \eqref{eq40}, we would have, in particular,
\begin{equation}\label{eq41}%{eq41}
m'(A)\leq m''(A^{\vep_0})+\vep_0,
\end{equation}
which becomes simply
\begin{equation}\label{eq42}%{eq42}
1=\fr{\text{no. of }a'_i\text{ in }\cub{a'_1,\dots,a'_{n'}}}{n'}\leq\fr{\text{no. of }a''_j\text{ in }\bigcup_{i=1}^{n'}\bD_{\vep_0}(a'_i)}{n''}+\vep_0
\end{equation}
when taking $A=\cub{a'_1,\dots,a'_{n'}}$ so that $A^{\vep_0}=\bigcup_{i=1}^{n'}\bD_{\vep_0}(a'_i)$. Now, if we suppose that less than $\floor{qn''}$ of $a''_1,\dots,a''_{n''}$ satisfy
\begin{equation}\label{eq43}%{eq43}
\min_{1\leq i\leq n'}\abs{a'_i-a''_j}<\vep_0,\text{ \ or, equivalently, \ }\abs{a'_i-a''_j}<\vep_0\text{ \ for some \ }i,
\end{equation}
then \eqref{eq42} would proceed as
$$1=\cdots\leq\cdots<\fr{\floor{qn''}}{n''}+\vep_0\leq q+\vep_0,$$
which contradicts \eqref{eq40}. Hence, we conclude that at least $\floor{qn''}$ of $a''_1,\dots,a''_{n''}$ satisfy \eqref{eq43} under assumption \eqref{eq40}, and then consequently:

%---------- COROLLARY 8 ----------
\bigskip
\begin{cor}\label{cor8}
Let $0<q<1$ {\color{b}and $0<\vep_0<1-q$ be fixed. Let $Z_1,\dots,Z_n$ and $W_{n,1}$, ..., $W_{n,n-k}$ be as in \eqref{eq1} and \eqref{eq3} respectively. Consider the event $E_n$ consisting of all the samples $\ome\in\Ome$ that satisfies:
\begin{quote}
At least $\floor{qn}$ of $Z_1(\ome),\dots,Z_n(\ome)$ satisfy
$$\min_{1\leq i\leq n-k}\abs{W_{n,i}(\ome)-Z_j(\ome)}<\vep_0,$$
and at least $\floor{q(n-k)}$ of $W_{n,1}(\ome),\dots,W_{n,n-k}(\ome)$ satisfy
$$\min_{1\leq j\leq n}\abs{Z_j(\ome)-W_{n,i}(\ome)}<\vep_0.$$
\end{quote}
If the weak convergence \eqref{eq6} holds almost surely, then $\bP(E_n\text{, ev.})=1$.} In particular, Theorem \ref{thm1} (which is the sum of Corollaries \ref{cor4} and \ref{cor6}) and Corollary \ref{cor5} are such cases.
\end{cor}

%---------- PROOF OF COROLLARY 8 ----------
\bigskip
\begin{proof}
Consider the empirical measure
$$\mu^{(0)}_{S,n}(\ome):=\nth{n}\sum_{j=1}^n\delta_{Z_j(\ome)}$$
of $Z_1(\ome),\dots,Z_n(\ome)$. For any bounded continuous function $f:\bC\to\bR$, since $\mu_S\sim Z=Z_1,Z_2,\dots$ are i.i.d., their post-compositions with $f$,
$$f\circ Z=f\circ Z_1,f\circ Z_2,\dots\phantom{}:(\Ome,\bP)\to\bR,$$
are also i.i.d., so by Kolmogorov's strong law of large numbers,
\begin{align*}
\int_{\ze\in\bC}f(\ze)\,d\mu^{(0)}_{S,n}(\ze)&=\nth{n}\sum_{j=1}^{n}\int_{\ze\in\bC}f(\ze)\,d\delta_{Z_j}(\ze)=\nth{n}\sum_{j=1}^{n}f\circ Z_j\\
&\to\bE\sqb{f\circ Z}=\int_{\ze\in\bC}f(\ze)\,d\mu_S(\ze)\antti\text{ \ a.s..}
\end{align*}
Therefore,
\begin{equation}\label{eq45}%{eq45}
\mu^{(0)}_{S,n}\wto\mu_S\antti\text{ \ a.s..}
\end{equation}
{\color{b}By that \eqref{eq6} holds almost surely, we have
\begin{equation}\label{eq105}%{eq105}
\pi\rob{\mu_{S,n},\mu_S}<\fr{\vep_0}{2}\text{ \ for sufficiently large \ }n\text{ \ a.s..}
\end{equation}
Combining \eqref{eq45} and \eqref{eq105}, we would have
$$\pi\rob{\mu_{S,n},\mu^{(0)}_{S,n}}<\vep_0\text{ \ for sufficiently large \ }n\text{ \ a.s..}$$
The required result then simply follows from using the preceding discussion with
$$m'=\mu_{S,n}\text{ \ and \ }m''=\mu^{(0)}_{S,n},\quad\text{and}\quad m'=\mu^{(0)}_{S,n}\text{ \ and \ }m''=\mu_{S,n}$$
in \eqref{eq40} respectively.}
\end{proof}

%---------- FIGURE 1 ----------
\bigskip
\begin{figure}[h]
\centering
\includegraphics[width=0.5\textwidth]{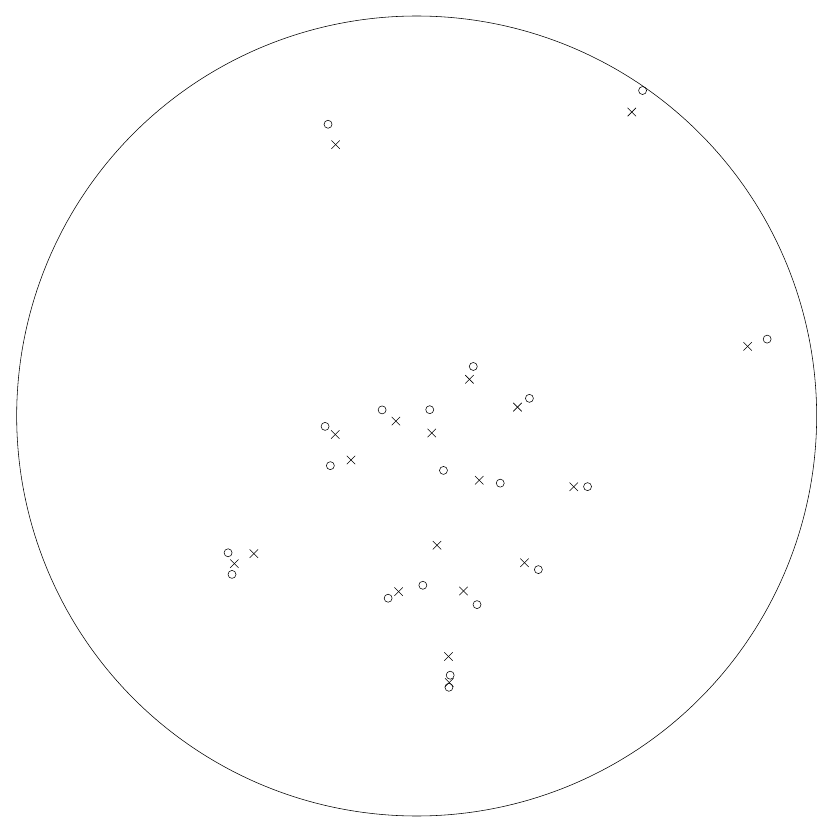}\includegraphics[width=0.5\textwidth]{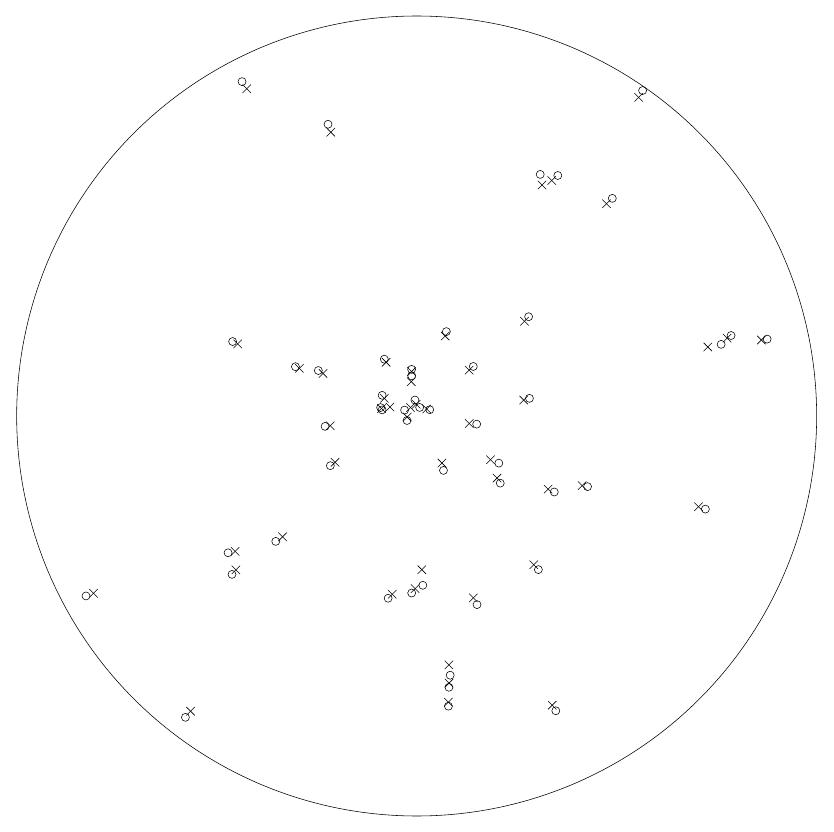}
\caption{The zeros $Z_1,\dots,Z_n\in S=\ol\bD$ (circle) of $P_n$ and its critical points {\color{b}$W^{(1)}_{n,1},\dots,W^{(1)}_{n,n-1}$} (cross) ($n=20,50$ resp.)\label{fig1}}%{fig1}
\end{figure}

%---------- REMARK 3 ----------
\bigskip
\begin{remark}
\begin{enumerate}
\item[(i)]	A remark similar to Corollary \ref{cor8} {\color{b}but regarding that the weak convergence \eqref{eq6} holds in probability has appeared} in \cite{Ka} without proof.
\item[(ii)]	Figure \ref{fig1} illustrates an instance that when more and more zeros are added to a polynomial (i.e. as $n\to\infty$), a majority of the critical points tend to get closer and closer to the zeros. Corollary \ref{cor8} assures high chances of similar phenomena for the zeros of higher order, polar and Sz.-Nagy's generalized derivatives, instead of the critical points. These somehow respond probabilistically to high degree situation of the old open problem Sendov's conjecture in geometry of polynomials, and suggest extensions to other types of derivatives.
\end{enumerate}
\end{remark}

% =======================================================================

\end{document}